\documentclass[12pt]{article}
\usepackage{latexsym,amssymb,enumerate,float,geometry,cite}
\newtheorem{theorem}{Theorem}
\newtheorem{corollary}[theorem]{Corollary}
\newtheorem{conjecture}[theorem]{Conjecture}
\newtheorem{claim}{Claim}

\usepackage{setspace}\onehalfspace

\newcommand{\smallqed}{{\tiny ($\Box$)}}
\newcommand{\cart}{\, \Box \,}

\newenvironment{unnumbered}[1]{\trivlist \item [\hskip \labelsep {\bf
#1}]\ignorespaces\it}{\endtrivlist}

\begin{document}

\title{Hereditary Equality of Domination\\ and Exponential Domination}
\author{Michael A. Henning$^1$, Simon J\"{a}ger$^2$, Dieter Rautenbach$^2$}
\date{}
\maketitle
\vspace{-10mm}
\begin{center}
{\small
$^1$ Department of Pure and Applied Mathematics, University of Johannesburg,\\
Auckland Park, 2006, South Africa, \texttt{mahenning@uj.ac.za}\\[3mm]
$^2$
Institute of Optimization and Operations Research, Ulm University,\\
Ulm, Germany, \texttt{simon.jaeger,dieter.rautenbach@uni-ulm.de}}
\end{center}

\begin{abstract}
We characterize a large subclass of the class of those graphs $G$
for which the exponential domination number of $H$ equals the domination number of $H$ for every induced subgraph $H$ of $G$.
\end{abstract}

\section{Introduction}

Domination in graphs is an important area within graph theory,
and an astounding variety of different domination parameters are known \cite{hhs}. 
Essentially all of these parameters involve merely local conditions, 
which makes them amenable to similar approaches and arguments.
In \cite{ddems} Dankelmann et al.~introduce a truly non-local variant of domination, the so-called exponential domination, 
where the influence of vertices extends to any arbitrary distance within the graph but decays exponentially with that distance.
There is relatively few research concerning exponential domination \cite{abcvy,aa,bor1,bor2}, and even apparently basic results require new and careful arguments.

Bessy et al.~\cite{bor2} show that the exponential domination number is APX-hard for subcubic graphs and describe an efficient algorithm for subcubic trees, but the complexity for general trees is unknown.
It is not even known how to decide efficiently 
for a given tree $T$
whether its exponential domination number $\gamma_e(T)$
equals its domination number $\gamma(T)$.
In \cite{hjr} we study relations between the different parameters of exponential domination and domination. Next to several bounds, we obtain a constructive characterization of the subcubic trees $T$ with $\gamma_e(T)=\gamma(T)$.
In view of the efficient algorithms to determine both parameters for such trees, the existence of a constructive characterization is not surprising, but, as said a few lines above, already for general trees all our techniques completely fail.

The difficulty to decide whether $\gamma_e(G)=\gamma(G)$ for a given graph $G$ motivates the study of the hereditary class ${\cal G}$ of graphs that satisfy this equality, that is, 
${\cal G}$ is the set of those graphs $G$ such that $\gamma_e(H)=\gamma(H)$ for every induced subgraph $H$ of $G$.
As for the well-known perfect graphs, the class ${\cal G}$ can be characterized by minimal forbidden induced subgraphs.

In the present paper we obtain such a characterization for a large subclass of ${\cal G}$,
and pose several related conjectures.

\medskip

\noindent Before we proceed to our results, we collect some notation.
We consider finite, simple, and undirected graphs, and use standard terminology.
The vertex set and the edge set of a graph $G$ are denoted by $V(G)$ and $E(G)$, respectively.
The order $n(G)$ of $G$ is the number of vertices of $G$.
For a vertex $u$ of $G$,
the neighborhood of $u$ in $G$
and the degree of $u$ in $G$ are denoted by 
$N_G(u)$ and $d_G(u)$, respectively.
The distance ${\rm dist}_G(X,Y)$ between two sets $X$ and $Y$ of vertices in $G$ is the minimum length of a path in $G$ 
between a vertex in $X$ and a vertex in $Y$.
If no such path exists, then let ${\rm dist}_G(X,Y)=\infty$.

Let $D$ be a set of vertices of a graph $G$.
The set $D$ is a dominating set of $G$ \cite{hhs}
if every vertex of $G$ not in $D$ has a neighbor in $D$.
The domination number $\gamma(G)$ of $G$ is the minimum order of a dominating set of $G$.
For two vertices $u$ and $v$ of $G$,
let ${\rm dist}_{(G,D)}(u,v)$ be the minimum length of a path $P$ in $G$ between $u$ and $v$ 
such that $D$ contains exactly one endvertex of $P$ but no internal vertex of $P$.
If no such path exists, then let ${\rm dist}_{(G,D)}(u,v)=\infty$.
Note that, if $u$ and $v$ are distinct vertices in $D$, 
then ${\rm dist}_{(G,D)}(u,u)=0$ and ${\rm dist}_{(G,D)}(u,v)=\infty$.
For a vertex $u$ of $G$, let 
\begin{eqnarray*}\label{ew}
w_{(G,D)}(u)=\sum\limits_{v\in D}\left(\frac{1}{2}\right)^{{\rm dist}_{(G,D)}(u,v)-1},
\end{eqnarray*}
where $\left(\frac{1}{2}\right)^{\infty}=0$.
Dankelmann et al.~\cite{ddems} define 
the set $D$ to be an {\it exponential dominating set} of $G$ if $w_{(G,D)}(u)\geq 1$ for every vertex $u$ of $G$, 
and the {\it exponential domination number $\gamma_e(G)$} of $G$ 
as the minimum order of an exponential dominating set of $G$.
Note that $w_{(G,D)}(u)=2$ for $u\in D$,
and that $w_{(G,D)}(u)\geq 1$ for every vertex $u$ that has a neighbor in $D$,
which implies $\gamma_e(G)\leq \gamma(G)$.

\medskip

\noindent The following Figure  \ref{figure1} contains forbidden induced subgraphs that relate to the considered subclasses of ${\cal G}$.

\begin{figure}[H]
\begin{center}
$\mbox{}$\hfill
\unitlength 1mm 
\linethickness{0.4pt}
\ifx\plotpoint\undefined\newsavebox{\plotpoint}\fi 
\begin{picture}(16,14)(0,0)
\put(5,13){\circle*{2}}
\put(15,13){\circle*{2}}
\put(5,13){\line(1,0){10}}
\put(10,0){\makebox(0,0)[cc]{$K_3$}}
\put(10,5){\circle*{2}}
\multiput(15,13)(-.033557047,-.053691275){149}{\line(0,-1){.053691275}}
\multiput(10,5)(-.033557047,.053691275){149}{\line(0,1){.053691275}}
\end{picture}
\hfill
\linethickness{0.4pt}
\ifx\plotpoint\undefined\newsavebox{\plotpoint}\fi 
\begin{picture}(16,24)(0,0)
\put(5,13){\circle*{2}}
\put(15,13){\circle*{2}}
\put(10,0){\makebox(0,0)[cc]{$K_{2,3}$}}
\put(10,5){\circle*{2}}
\multiput(15,13)(-.033557047,-.053691275){149}{\line(0,-1){.053691275}}
\multiput(10,5)(-.033557047,.053691275){149}{\line(0,1){.053691275}}
\put(15,23){\circle*{2}}
\put(5,23){\circle*{2}}
\put(15,23){\line(0,-1){10}}
\put(5,23){\line(0,-1){10}}
\put(15,13){\line(-1,1){10}}
\put(15,23){\line(-1,-1){10}}
\end{picture}
\hfill
\linethickness{0.4pt}
\ifx\plotpoint\undefined\newsavebox{\plotpoint}\fi 
\begin{picture}(21,16)(0,0)
\put(0,5){\circle*{2}}
\put(10,5){\circle*{2}}
\put(20,5){\circle*{2}}
\put(0,15){\circle*{2}}
\put(10,15){\circle*{2}}
\put(20,15){\circle*{2}}
\put(10,15){\line(-1,0){10}}
\put(20,15){\line(-1,0){10}}
\put(0,15){\line(0,-1){10}}
\put(10,15){\line(0,-1){10}}
\put(20,15){\line(0,-1){10}}
\put(0,5){\line(1,0){10}}
\put(10,5){\line(1,0){10}}
\put(10,0){\makebox(0,0)[cc]{$P_2 \cart P_3$}}
\end{picture}
\hfill
\linethickness{0.4pt}
\ifx\plotpoint\undefined\newsavebox{\plotpoint}\fi 
\begin{picture}(16,24)(0,0)
\put(5,13){\circle*{2}}
\put(15,13){\circle*{2}}
\put(5,13){\line(1,0){10}}
\put(10,0){\makebox(0,0)[cc]{$B$}}
\put(10,5){\circle*{2}}
\multiput(15,13)(-.033557047,-.053691275){149}{\line(0,-1){.053691275}}
\multiput(10,5)(-.033557047,.053691275){149}{\line(0,1){.053691275}}
\put(15,23){\circle*{2}}
\put(5,23){\circle*{2}}
\put(15,23){\line(0,-1){10}}
\put(5,23){\line(0,-1){10}}
\end{picture}
\hfill
\linethickness{0.4pt}
\ifx\plotpoint\undefined\newsavebox{\plotpoint}\fi 
\begin{picture}(16,22)(0,0)
\put(5,13){\circle*{2}}
\put(15,13){\circle*{2}}
\put(5,13){\line(1,0){10}}
\put(10,0){\makebox(0,0)[cc]{$D$}}
\put(10,5){\circle*{2}}
\multiput(15,13)(-.033557047,-.053691275){149}{\line(0,-1){.053691275}}
\multiput(10,5)(-.033557047,.053691275){149}{\line(0,1){.053691275}}
\put(10,21){\circle*{2}}
\multiput(15,13)(-.033557047,.053691275){149}{\line(0,1){.053691275}}
\multiput(10,21)(-.033557047,-.053691275){149}{\line(0,-1){.053691275}}
\end{picture}
\hfill
\linethickness{0.4pt}
\ifx\plotpoint\undefined\newsavebox{\plotpoint}\fi 
\begin{picture}(16,22)(0,0)
\put(5,13){\circle*{2}}
\put(15,13){\circle*{2}}
\put(5,13){\line(1,0){10}}
\put(10,0){\makebox(0,0)[cc]{$K_4$}}
\put(10,5){\circle*{2}}
\multiput(15,13)(-.033557047,-.053691275){149}{\line(0,-1){.053691275}}
\multiput(10,5)(-.033557047,.053691275){149}{\line(0,1){.053691275}}
\put(10,21){\circle*{2}}
\multiput(15,13)(-.033557047,.053691275){149}{\line(0,1){.053691275}}
\multiput(10,21)(-.033557047,-.053691275){149}{\line(0,-1){.053691275}}
\put(10,21){\line(0,-1){16}}
\end{picture}
\hfill$\mbox{}$
\end{center}
\caption{The graphs $K_3$, $K_{2,3}$, $P_2 \cart P_3$, $B$ (bull), $D$ (diamond), and $K_4$.}\label{figure1}
\end{figure}
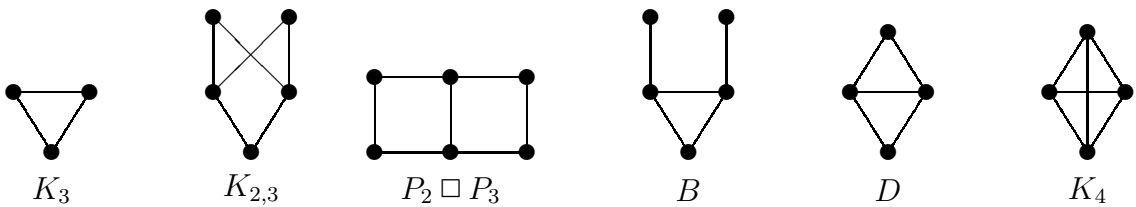
\noindent Our main result is the following, where $P_7$ and $C_7$ denote the path and the cycle of order $7$.

\begin{figure}[H]
\begin{center}
$\mbox{}$\hfill
\unitlength 1mm 
\linethickness{0.4pt}
\ifx\plotpoint\undefined\newsavebox{\plotpoint}\fi 
\begin{picture}(21,16)(0,0)
\put(0,5){\circle*{2}}
\put(10,5){\circle*{2}}
\put(20,5){\circle*{2}}
\put(0,15){\circle*{2}}
\put(10,15){\circle*{2}}
\put(20,15){\circle*{2}}
\put(10,15){\line(-1,0){10}}
\put(20,15){\line(-1,0){10}}
\put(0,15){\line(0,-1){10}}
\put(10,15){\line(0,-1){10}}
\put(20,15){\line(0,-1){10}}
\put(10,0){\makebox(0,0)[cc]{$F_1$}}
\end{picture}
\hfill
\linethickness{0.4pt}
\ifx\plotpoint\undefined\newsavebox{\plotpoint}\fi 
\begin{picture}(31,16)(0,0)
\put(0,5){\circle*{2}}
\put(10,5){\circle*{2}}
\put(0,15){\circle*{2}}
\put(10,15){\circle*{2}}
\put(20,15){\circle*{2}}
\put(10,15){\line(-1,0){10}}
\put(20,15){\line(-1,0){10}}
\put(0,15){\line(0,-1){10}}
\put(10,15){\line(0,-1){10}}
\put(15,0){\makebox(0,0)[cc]{$F_2$}}
\put(0,5){\line(1,0){10}}
\put(30,15){\circle*{2}}
\put(30,5){\circle*{2}}
\put(19,15){\line(1,0){11}}
\put(30,15){\line(0,-1){10}}
\end{picture}
\hfill
\linethickness{0.4pt}
\ifx\plotpoint\undefined\newsavebox{\plotpoint}\fi 
\begin{picture}(31,16)(0,0)
\put(0,5){\circle*{2}}
\put(10,5){\circle*{2}}
\put(0,15){\circle*{2}}
\put(10,15){\circle*{2}}
\put(20,15){\circle*{2}}
\put(10,15){\line(-1,0){10}}
\put(20,15){\line(-1,0){10}}
\put(0,15){\line(0,-1){10}}
\put(10,15){\line(0,-1){10}}
\put(15,0){\makebox(0,0)[cc]{$F_3$}}
\put(0,5){\line(1,0){10}}
\put(20,5){\circle*{2}}
\put(10,5){\line(1,0){10}}
\put(20,5){\line(2,1){10}}
\put(30,10){\line(-2,1){10}}
\put(30,10){\circle*{2}}
\end{picture}
\hfill$\mbox{}$\\[10mm]
$\mbox{}$\hfill
\linethickness{0.4pt}
\ifx\plotpoint\undefined\newsavebox{\plotpoint}\fi 
\begin{picture}(41,16)(0,0)
\put(0,10){\circle*{2}}
\put(10,5){\circle*{2}}
\put(10,15){\circle*{2}}
\put(20,10){\circle*{2}}
\put(30,15){\circle*{2}}
\put(30,5){\circle*{2}}
\put(40,10){\circle*{2}}
\put(40,10){\line(-2,1){10}}
\put(30,15){\line(-2,-1){10}}
\put(20,10){\line(-2,-1){10}}
\put(10,5){\line(-2,1){10}}
\put(0,10){\line(2,1){10}}
\put(10,15){\line(2,-1){10}}
\put(20,10){\line(2,-1){10}}
\put(30,5){\line(2,1){10}}
\put(20,0){\makebox(0,0)[cc]{$F_4$}}
\end{picture}
\hfill
\linethickness{0.4pt}
\ifx\plotpoint\undefined\newsavebox{\plotpoint}\fi 
\begin{picture}(41,16)(0,0)
\put(10,10){\circle*{2}}
\put(20,15){\circle*{2}}
\put(20,5){\circle*{2}}
\put(30,10){\circle*{2}}
\put(30,10){\line(-2,1){10}}
\put(20,15){\line(-2,-1){10}}
\put(10,10){\line(2,-1){10}}
\put(20,5){\line(2,1){10}}
\put(0,15){\circle*{2}}
\put(0,5){\circle*{2}}
\put(20,15){\line(-1,0){20}}
\put(0,15){\line(0,-1){10}}
\put(0,5){\line(1,0){20}}
\put(40,10){\circle*{2}}
\put(30,10){\line(1,0){10}}
\put(20,0){\makebox(0,0)[cc]{$F_5$}}
\end{picture}
\hfill$\mbox{}$
\end{center}
\caption{The graphs $P_2 \cart P_3, F_1,\ldots,F_5$.}\label{figure2}
\end{figure}
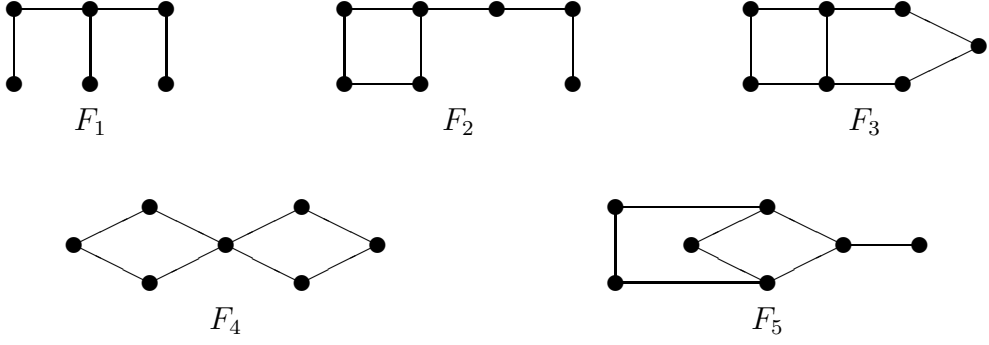

\begin{theorem}\label{theorem1}
If $G$ is a $\{B, D, K_4,K_{2,3},P_2 \cart P_3\}$-free graph,
then $\gamma(H) = \gamma_e(H)$ for every induced subgraph $H$ of $G$
if and only if $G$ is
$\{ P_7, C_7, F_1,\dots, F_5\}$-free (cf. Figure \ref{figure2}).
\end{theorem}

\noindent Since $K_3$ is an induced subgraph of the graphs $B$, $D$, and $K_4$, 
Theorem \ref{theorem1} has the following immediate corollary.

\begin{corollary}\label{corollary1}
If $G$ is a $\{ K_3,K_{2,3},P_2 \cart P_3\}$-free graph,
then $\gamma(H) = \gamma_e(H)$ for every induced subgraph $H$ of $G$
if and only if $G$ is
$\{ P_7, C_7, F_1,\dots, F_5\}$-free.
\end{corollary}
For the trees in ${\cal G}$, we achieve a complete characterization.

\begin{corollary}\label{corollary2}
If $T$ is a tree, 
then $\gamma(F) = \gamma_e(F)$ for every induced subgraph $F$ of $T$
if and only if $T$ is $\{ P_7,F_1\}$-free.
\end{corollary}
All proofs and our conjectures are postponed to the next section.

\section{Proofs and Conjectures}

In order to prove Theorem \ref{theorem1}, 
we first give a complete and independent proof of Corollary \ref{corollary1},
and then prove Theorem \ref{theorem1} relying on this partial result.

\medskip

\noindent {\it Proof of Corollary \ref{corollary1}:}
Since $\gamma(H)>\gamma_e(H)$ for every graph $H$ in $\{ P_7, C_7, F_1,\dots, F_5\}$,
the necessity follows. In order to prove the sufficiency, suppose that $G$ is a
$\{ K_3,K_{2,3},P_2 \cart P_3\}\cup\{ P_7, C_7, F_1,\dots, F_5\}$-free graph
with $\gamma(G)>\gamma_e(G)$ of minimum order.
By the choice of $G$, we have $\gamma(H) = \gamma_e(H)$ for every proper induced subgraph $H$ of $G$.
Clearly, $G$ is connected.
Since $\gamma_e(G)=1$ if and only if $\gamma(G)=1$,
we obtain $\gamma_e(G)\geq 2$ and $\gamma(G)\geq 3$.
Since $G$ is $\{ P_7,C_7\}$-free, the girth $g$ of $G$ is at most $6$,
where the girth of $G$ is the minimum length of a cycle of $G$.
Let $C:x_1x_2x_3\ldots x_gx_1$ be a shortest cycle of $G$,
where we consider the indices modulo $g$.
Let $R=V(G)\setminus V(C)$.

Suppose $g=6$.
Since $\gamma(C_6)=\gamma_e(C_6)=2$, some vertex $y$ in $R$ has a neighbor $x_i$ on $C$.
Since $g=6$, the vertex $y$ has no further neighbor on $C$, implying that $G[\{ y,x_{i-2},x_{i-1},x_i,x_{i+1},x_{i+2}\}] = F_1$, contradicting the fact that $G$ is $F_1$-free.
Hence, $g<6$.

Suppose $g=5$.
This implies that no vertex in $R$ has more than one neighbor on $C$.
If some vertex $z$ has distance $2$ from $V(C)$ in $G$ and $x_iyz$ is a path in $G$,
then $G[\{ z,y,x_{i-1},x_i,x_{i+1},x_{i+2}\}]=F_1$, which is a contradiction.
Hence, every vertex in $R$ has a unique neighbor on $C$. Suppose that there is some $i\in [5]$ such that
$x_i$ has a neighbor $y_i$ in $R$
and
$x_{i+1}$ has a neighbor $y_{i+1}$ in $R$. Since $g = 5$, we note that $y_i \ne y_{i+1}$ and that the vertex $y_i$ is not adjacent to $y_{i+1}$, implying that $G[\{ x_{i-2},x_{i-1},x_i,x_{i+1},y_i,y_{i+1}\}]=F_1$, which is a contradiction.
This implies the existence of some index $i\in [5]$ such that $\{ x_i,x_{i+2}\}$ is a dominating set of $G$,
which contradicts $\gamma(G) \geq 3$. Hence, $g \le 4$. Since $G$ is $K_3$-free, this implies that $g=4$.

Since $G$ is $\{K_3,K_{2,3}\}$-free, no vertex in $R$ has more than one neighbor on $C$, and
since $G$ is $F_2$-free, no vertex in $R$ has distance more than $2$ from $V(C)$.

Suppose that some vertex $z$ has distance $2$ from $V(C)$.
Let $x_1yz$ be a path in $G$.
Suppose that $x_2$ has a neighbor $u$ in $R$.
Recall that $u$ is not adjacent to a second vertex on $C$.
Since $G$ is $P_2 \cart P_3$-free, the vertex $u$ is not adjacent to $y$.
If $u$ is not adjacent to $z$, then $G[\{ u,x_1,x_2,x_4,y,z\}]=F_1$, which is a contradiction.
If $u$ is adjacent to $z$, then $G[V(C)\cup \{ u,y,z\}]=F_3$, which is a contradiction.
Hence, by symmetry, we obtain $d_G(x_2)=d_G(x_4)=2$.
Suppose that $x_1$ has a neighbor $u$ in $R\setminus \{ y\}$.
Since $G$ is $\{ K_3,K_{2,3}\}$-free, the vertex $u$ is not adjacent to any vertex in $\{ x_2,x_3,x_4,y\}$.
If $u$ is not adjacent to $z$, then $G[\{ x_1,x_2,x_3,u,y,z\}]=F_1$, which is a contradiction.
If $u$ is adjacent to $z$, then $G[V(C)\cup \{ u,y,z\}]=F_4$, which is a contradiction.
Hence, we obtain $d_G(x_1)=3$.

Since $\{ x_3,y\}$ is not a dominating set of $G$, and no vertex in $R$ has distance more than $2$ from $V(C)$,
the degrees of $x_1$, $x_2$, and $x_4$ imply the existence of a path $x_3uv$,
where $v$ has distance $2$ to $V(C)$,
and $v$ is not adjacent to $y$.
Since $G[\{ v,u,x_3,x_2,x_1,y,z\}]$ is neither $P_7$ nor $C_7$, the vertex $u$ is adjacent to $y$ or $z$.
If $u$ is adjacent to $z$, then, because $G$ is $K_3$-free,
$G[\{ u,v,y,z,x_3,x_2\}]=F_1$, which is a contradiction.
Hence, the vertex $u$ is adjacent to $y$.
If $v$ is adjacent to $z$, then
$G[\{ u,v,y,z,x_1,x_2,x_3\}]=F_3$, which is a contradiction.
Hence, the vertex $v$ is not adjacent to $z$, and
$G[\{ u,v,y,z,x_3,x_2\}]=F_1$, which is a contradiction.
Hence, every vertex in $R$ has a unique neighbor on $C$.

Since $\gamma(G)>2$, we may assume that $x_i$ has a neighbor $y_i$ in $R$ for $i\in [3]$.
Since $G$ is $P_2 \cart P_3$-free, the vertex $y_2$ is not adjacent to $y_1$ or $y_3$.
Since $G$ is $F_5$-free, the vertex $y_1$ is not adjacent to $y_3$.
Now, $G[\{ x_1,x_2,x_3,y_1,y_2,y_3\}]=F_1$, which is a contradiction,
and completes the proof. $\Box$

\medskip

\noindent With Corollary \ref{corollary1} at hand, we now proceed to the proof of Theorem \ref{theorem1}.

\medskip

\noindent {\it Proof of Theorem \ref{theorem1}:}
The necessity follows as above.
In order to prove the sufficiency, suppose that $G$ is a
$\{ B,D, K_4,K_{2,3},P_2 \cart  P_3\}\cup\{ P_7, C_7, F_1,\dots, F_5\}$-free graph
with $\gamma(G)>\gamma_e(G)$ of minimum order.
By the choice of $G$, we have $\gamma(H) = \gamma_e(H)$ for every proper induced subgraph $H$ of $G$.
Clearly, $G$ is connected.
Since $\gamma_e(G)=1$ if and only if $\gamma(G)=1$,
we obtain $\gamma_e(G)\geq 2$ and $\gamma(G)\geq 3$.

By Corollary \ref{corollary1}, $G$ is not $K_3$-free, that is, the girth of $G$ is $3$.

We proceed with a series of claims.
Let $F$ be the graph that is obtained from the triangle $x_1x_2x_3$ and the path $y_1y_2y_3$ by adding the edge $x_1y_1$.

\begin{claim}\label{claim1}
$F$ is not an induced subgraph of $G$.
\end{claim}
{\it Proof of Claim~\ref{claim1}:}
Suppose that $F$ is an induced subgraph of $G$.
Since $\{x_1, y_2\}$ is not a dominating set of $G$, there is a vertex $u$ at distance~$2$ from the set $\{x_1,y_2\}$ in $G$. 

We proceed with three subclaims.

\begin{unnumbered}{Claim~\ref{claim1}.1}
The vertex $u$ is not adjacent to $x_2$ or $x_3$.
\end{unnumbered}
{\it Proof of Claim~\ref{claim1}.1.} Suppose that $u$ is adjacent to $x_2$.
Since $G$ is $D$-free, the vertex $u$ is not adjacent to $x_3$, and,
since $G$ is $B$-free, $u$ is adjacent to $y_1$.
If $u$ is not adjacent to $y_3$, then $G[\{u,x_1, x_3, y_1, y_2, y_3\}]=F_1$, which is a contradiction.
If $u$ is adjacent to $y_3$, then $G[\{u,x_1, x_2, y_1, y_2, y_3\}]=P_2 \cart P_3$, which is a contradiction.
Hence, by symmetry, we obtain that $u$ is not adjacent to $x_2$ or $x_3$.~\smallqed

\begin{unnumbered}{Claim~\ref{claim1}.2}
The vertex $u$ is not adjacent to $y_1$.
\end{unnumbered}
{\it Proof of Claim~\ref{claim1}.2.} Suppose that $u$ is adjacent to $y_1$.
Since $G$ is $F_1$-free, the vertex $u$ is adjacent to $y_3$.
Since $\{x_1, y_3\}$ is not a dominating set of $G$,
there is a vertex $v$ at distance~$2$ from the set $\{x_1,y_3\}$.
Suppose that $v$ is adjacent to $x_2$.
Since $G$ is $D$-free, the vertex $v$ is not adjacent to $x_3$, and,
since $G$ is $B$-free, $v$ is adjacent to $y_1$.
If $v$ is adjacent to $u$, then $G[\{u,v,x_2, y_1, y_3\}]=B$, which is a contradiction.
Hence, $v$ is not adjacent to $u$, and, by symmetry, $v$ is also not adjacent to $y_2$.  
Therefore, $G[\{u,v,x_1, x_2, y_1, y_2, y_3\}]=F_4$, which is a contradiction.
Thus, by symmetry, $v$ is not adjacent to $x_2$ or $x_3$.
Next, suppose that $v$ is adjacent to $y_1$.
Since $G$ is $F_1$-free, the vertex $v$ is adjacent to both $u$ and $y_2$, which yields the contradiction $G[\{u,v,y_1, y_2\}]=D$.
Thus, $v$ is not adjacent to $y_1$.
Suppose that $v$ is adjacent to $u$.
If $v$ is adjacent to $y_2$, then $G[\{u,v, y_1,y_2, y_3\}]=K_{2,3}$, which is a contradiction.
If $v$ is not adjacent to $y_2$, then $G[\{u,v, x_1, x_2,y_1,y_2\}]=F_1$, which is a contradiction.
Thus, by symmetry, $v$ is not adjacent to $u$ or $y_2$, implying that $v$ is at distance~$2$ from the set $\{ u,x_1,x_2,x_3,y_1,y_2,y_3\}$.

Since the vertex $v$ is at distance~$2$ from the set $\{x_1,y_3\}$ in $G$, 
there is a neighbor $v'$ of $v$, that is adjacent to $x_1$ or to $y_3$ or to both $x_1$ and $y_3$.
First, suppose that $v'$ is not adjacent to $x_1$, implying that $v'$ is adjacent to $y_3$.
Suppose that $v'$ is adjacent to $y_1$. If $v'$ is adjacent to $u$, then $G[\{u,v',y_1,y_3\}] = D$, which is a contradiction. Thus, by symmetry, $v'$ is adjacent to neither $u$ nor $y_2$, implying that
$G[\{u,v', y_1, y_2, y_3\}]=K_{2,3}$, which is a contradiction.
Thus, $v'$ is not adjacent to $y_1$.
Since $G$ is $B$-free, $v'$ is not adjacent to $u$ or $y_2$, which yields the contradiction $G[\{v,v',x_1, x_2, y_1, y_2, y_3\}]=P_7$.
Therefore, $v'$ is adjacent to $x_1$.

Since $G$ is $\{B,K_4\}$-free, the vertex $v'$ is not adjacent to $x_2$ or $x_3$.
Since $G$ is $B$-free, the vertex $v'$ is not adjacent to $y_1$.
Suppose that $v'$ is not adjacent to $y_3$.
If $v'$ is not adjacent to $u$, then $G[\{u,v, v', x_1, x_2, y_1\}]=F_1$, which is a contradiction.
Thus, by symmetry, $v'$ is adjacent to both $u$ and $y_2$, which yields the contradiction $G[\{u,v', y_1, y_2, y_3\}]=K_{2,3}$.
Thus, $v'$ is adjacent to $y_3$.
Since $G$ is $F_1$-free, the vertex $v'$ is adjacent to both $u$ and $y_2$, implying that $G[\{u,v,v', y_1, y_3\}]=B$, which is a contradiction. Therefore, $u$ is not adjacent to $y_1$.~\smallqed

\begin{unnumbered}{Claim~\ref{claim1}.3}
The vertex $u$ is not adjacent to $y_3$.
\end{unnumbered}
{\it Proof of Claim~\ref{claim1}.3.} Suppose that $u$ is adjacent to $y_3$. Since $\{x_1, y_3\}$ is not a dominating set of $G$, 
there is a vertex $v$ at distance~$2$ from the set $\{x_1,y_3\}$ in $G$.

Suppose that $v$ is adjacent to $x_2$.
Since $G$ is $D$-free, the vertex $v$ is not adjacent to $x_3$.
Since $G$ is $B$-free, $v$ is adjacent to $y_1$.
If $v$ is not adjacent to $y_2$, then $G[\{v,x_1,x_3,y_1, y_2, y_3\}]=F_1$, which is a contradiction.
If $v$ is adjacent to $y_2$, then we get the contradiction $G[\{x_2,v,y_1, y_2, y_3\}]=B$.
Therefore, by symmetry, $v$ is not adjacent to $x_2$ or $x_3$.

Next, suppose that $v$ is adjacent to $y_1$.
If $v$ is not adjacent to $y_2$, then $G[\{v,x_1,x_2,y_1, y_2, y_3\}]=F_1$, which is a contradiction.
If $v$ is adjacent to $y_2$, then $G[\{v,x_1,y_1, y_2, y_3\}]=B$, which is a contradiction.
Thus, $v$ is not adjacent to $y_1$.

Next, suppose that $v$ is adjacent to $y_2$.
If $v$ is not adjacent to $u$, then $G[\{u,v,x_1,y_1, y_2, y_3\}]=F_1$, which is a contradiction.
If $v$ is adjacent to $u$, then $G[\{u,v,x_1,x_2,y_1, y_2, y_3\}]=F_2$, which is a contradiction.
Therefore, $v$ is not adjacent to $y_2$, implying that $v$ is at distance~$2$ from the set $\{x_1,x_2,x_3,y_1,y_2,y_3\}$.

If $v$ is adjacent to $u$, then $G[\{u,v, x_1,x_2,y_1, y_2, y_3\}]=P_7$, which is a contradiction. Hence, $v$ is not adjacent to $u$.
Since the vertex $v$ is at distance~$2$ from the set $\{x_1,y_3\}$ in $G$, 
there is a neighbor $v'$ of $v$ that is adjacent to $x_1$ or to $y_3$ or to both $x_1$ and $y_3$. 
Note that $v' \ne u$.

First, suppose that $v'$ is not adjacent to $x_1$, implying that $v'$ is adjacent to $y_3$. If $v'$ is not adjacent to $y_2$, then analogous arguments as in Claim~1.1 and Claim~1.2 (with the vertex $u$ replaced by the vertex $v'$) show that $y_3$ is the only vertex in the set $\{x_1,x_2,x_3,y_1,y_2,y_3\}$ that is adjacent to $v'$. This in turn implies that $G[\{v,v', x_1,x_2,y_1, y_2, y_3\}]=P_7$, which is a contradiction. Hence, $v'$ is adjacent to $y_2$. If $v'$ is adjacent to $y_1$, then $G[\{v',y_1,y_2,y_3\}] = D$, which is a contradiction. Thus, $v'$ is not adjacent to $y_1$, implying that
$G[\{v,v', y_1, y_2, y_3\}]=B$, which is a contradiction. Therefore, $v'$ is adjacent to $x_1$.

Since $G$ is $\{D,K_4\}$-free, the vertex $v'$ is not adjacent to $x_2$ or $x_3$.
If $v'$ is adjacent to $y_1$, then $G[\{v,v',x_1,x_2,y_1\}] = B$, which is a contradiction. Thus, $v'$ is not adjacent to $y_1$.
If $v'$ is not adjacent to $y_2$, then $G[\{v,v',x_1,x_2,y_1,y_2\}] = F_1$, which is a contradiction. Thus, $v'$ is adjacent to $y_2$.
If $v'$ is not adjacent to $y_3$, then $G[\{v, v', x_1, x_2, y_2,y_3\}]=F_1$, which is a contradiction. Thus, $v'$ is adjacent to $y_3$, implying that $G[\{v,v',y_1,y_2,y_3\}] = B$, which is a contradiction. Therefore, $u$ is not adjacent to $y_3$.~\smallqed

\medskip

\noindent We return to the proof of Claim~\ref{claim1}. By Claim~\ref{claim1}.1, Claim~\ref{claim1}.2 and Claim~\ref{claim1}.3, the vertex $u$ is at distance~$2$ from the set $\{x_1,x_2,x_3,y_1,y_2,y_3\}$. Since the vertex $u$ is at distance~$2$ from the set $\{x_1,y_2\}$ in $G$, there is a neighbor $u'$ of $u$ that is adjacent to $x_1$ or to $y_2$ or to both $x_1$ and $y_2$.
First, suppose that $u'$ is adjacent to $x_1$. 
Analogously as above, since $G$ is $\{B,D,K_4\}$-free, the vertex $u'$ is not adjacent to $x_2$, $x_3$ and $y_1$. 
If $u'$ is not adjacent to $y_2$, then $G[\{u,u',x_1,x_2,y_1, y_2\}]=F_1$, which is a contradiction. 
Thus, $u'$ is adjacent to $y_2$. If $u'$ is adjacent to $y_3$, then $G[\{u,u',y_1, y_2,y_3\}]=B$, while, if $u'$ is not adjacent to $y_3$, then $G[\{u,u',x_1,x_2,y_2,y_3\}]=F_1$. Since both cases produce a contradiction, we deduce that $u'$ is not adjacent to $x_1$, implying that $u'$ is adjacent to $y_2$. Since $G$ is $B$-free, $u'$ is not adjacent to $y_1$.
If $u'$ is not adjacent to $y_3$, then $G[\{u, u',x_1,y_1,y_2,y_3\}]=F_1$, which is a contradiction.
If $u'$ is adjacent to $y_3$, then $G[\{u, u', y_1,y_2,y_3\}]=B$, which is a contradiction.
This completes the proof of Claim \ref{claim1}.~$\Box$

\begin{claim}\label{claim2}
If $C$ is an arbitrary triangle in $G$, then every vertex is within distance $2$ from $V(C)$.
\end{claim}
{\it Proof of Claim~\ref{claim2}:}
Let $C \colon x_1x_2x_3$ be a triangle in $G$.
Suppose that there is a vertex $y_3$ at distance~$3$ from $V(C)$ in $G$. Let $x_1y_1y_2y_3$ be a shortest path in $G$ from $y_3$ to $V(C)$. Since $G$ is $\{D,K_4\}$-free, the vertex $y_1$ is adjacent to neither $x_2$ nor $x_3$, implying that $F$ is an induced subgraph of $G$, which contradicts Claim \ref{claim1}.~$\Box$

\begin{claim}\label{claim3}
Every triangle contains at least one vertex of degree exactly $2$ in $G$.
\end{claim}
{\it Proof of Claim~\ref{claim3}:}
Let $C \colon x_1x_2x_3$ be a triangle in $G$.
Suppose that every vertex on $C$ has degree at least $3$ in $G$.
Let $y_1, y_2, y_3\in V(G)\setminus V(C)$ be neighbors of $x_1, x_2, x_3$, respectively.
Since $G$ is $\{ D,K_4\}$-free, $x_i$ is the only neighbor of $y_i$ in $V(C)$ for $i\in [3]$.
Since $G$ is $B$-free, the vertices $y_1$, $y_2$ and $y_3$ induce a triangle $C'$ in $G$.
Suppose that there is a vertex $y\in V(G)\setminus (V(C)\cup V(C'))$ that is adjacent to a vertex on $C$, say $x_1$.
Since $G$ is $\{ D,K_4\}$-free, $x_1$ is the only neighbor of $y$ on $C$,
and $y$ is non-adjacent to some vertex $y_j$ on $C'$ with $j\in \{ 2,3\}$,
which implies the contradiction that $G[\{ x_1,x_2,x_3,y,y_j\}]=B$.
Hence, each vertex on $C$ has degree exactly~$3$ in $G$. 
By symmetry, each vertex on $C'$ has degree exactly~$3$ in $G$. 
Thus, $G = P_2 \cart C_3$, implying that $\gamma(G)=\gamma_e(G)=2$, which is a contradiction. 
This completes the proof of Claim \ref{claim3}.
$\Box$

\begin{claim}\label{claim4}
Every triangle contains two vertices of degree exactly $2$ in $G$.
\end{claim}
{\it Proof of Claim~\ref{claim4}:}
Let $C \colon x_1x_2x_3$ be a triangle in $G$ and let $R=V(G)\setminus V(C)$.
By Claim~\ref{claim3}, the triangle $C$ contains at least one vertex of degree exactly~$2$ in $G$. Renaming vertices if necessary, we may assume that $x_1$ has degree~$2$ in $G$. Suppose that both $x_2$ and $x_3$ have degree at least~$3$ in $G$.
Since $G$ is $D$-free, the vertices $x_2$ and $x_3$ have no common neighbor in $R$.
Further, since $G$ is $B$-free, every neighbor of $x_2$ in $R$ is adjacent to every neighbor of $x_3$ in $R$. Hence, since $G$ is $\{D,K_{2,3}\}$-free, the degrees of $x_2$ and $x_3$ are exactly~$3$ in $G$. 
Let $y_2$ and $y_3$ in $R$ be neighbors of $x_2$ and $x_3$, respectively.
Recall that $\gamma(G) \ge 3$. Let $w_2$ be a vertex not dominated by $\{x_2,y_3\}$, and let $w_3$ be a vertex not dominated by $\{x_3,y_2\}$. By Claim~\ref{claim2}, the vertex $w_2$ is within distance~$2$ from $V(C)$, implying that $w_2$ is adjacent to $y_2$. Analogously, the vertex $w_3$ is adjacent to $y_3$. Note that $w_2 \ne w_3$. If $w_2$ is adjacent to $w_3$, then
$G[\{w_2, w_3,x_2,x_3,y_2,y_3\}]=P_2 \cart P_3$. If $w_2$ is not adjacent to $w_3$, then $G[\{w_2, w_3,x_1,x_2,y_2,y_3\}]=F_1$.
Both cases produce a contradiction, which completes the proof of Claim~\ref{claim4}.
$\Box$

\medskip

\noindent Let $C \colon x_1x_2x_3$ be a triangle in $G$. 
By Claim \ref{claim4}, we may assume, renaming vertices if necessary,  that $x_2$ and $x_3$ have degree~$2$ in $G$.
 Since $\gamma(G) \ge 3$, the vertex $x_1$ does not dominate $V(G)$. Let $D_2 = V(G) \setminus N_G[x_1]$. 
 Claim~\ref{claim2} implies that every vertex in $D_2$ is at distance exactly~$2$ from $x_1$ in $G$. 
 Let $D_1$ be the set of neighbors in $V(G) \setminus D_2$ of the vertices in $D_2$. 
 Note that $D_1 \subset N_G(x_1)$. By Claim~\ref{claim4}, the set $D_1$ is independent. 

\begin{claim}\label{claim5}
Every vertex in $D_2$ has exactly one neighbor in $D_1$.
\end{claim}
{\it Proof of Claim~\ref{claim5}:} 
Since $D_1$ is an independent set, and, since $G$ is $K_{2,3}$-free, every vertex in $D_2$ has at most two neighbors in $D_1$. 
Suppose that a vertex $w_1$ in $D_2$ has two neighbors $y_1, y_2$ in $D_1$. 
Since $\{x_1, y_1\}$ is not a dominating set of $G$, there is a vertex $w_2\in D_2$ that is not adjacent to $y_1$.

\begin{unnumbered}{Claim~\ref{claim5}.1}
The vertex $w_2$ is not adjacent to $y_2$.
\end{unnumbered}
{\it Proof of Claim~\ref{claim5}.1.} Suppose that $w_2$ is adjacent to $y_2$.
Since $\{x_1, y_2\}$ is not a dominating set, 
there is a vertex $w_3$ in $D_2$ that is not adjacent to $y_2$. Suppose that $w_3$ is adjacent to $y_1$. If $w_3$ is not adjacent to $w_2$, then $G[\{x_1,x_2,y_1,y_2,w_2,w_3\}]=F_1$, which is a contradiction. Hence, $w_3$ is adjacent to $w_2$. If $w_3$ is adjacent to $w_1$, 
then, since $G$ is $D$-free, $w_1$ is not adjacent to $w_2$, implying that $G[\{x_1,y_1,w_1,w_2,w_3\}]=B$, which is a contradiction. Thus, $w_3$ is not adjacent to $w_1$. If $w_1$ is not adjacent to $w_2$, then  $G[\{x_1,x_2,y_1,w_1,w_2,w_3\}]=F_1$, while, if $w_1$ is adjacent to $w_2$, then  $G[\{x_1,y_2,w_1,w_2,w_3\}]=B$. Since both cases produce a contradiction, we deduce that $w_3$ is not adjacent to $y_1$. Since $G$ is $P_2 \cart P_3$-free, the vertex $w_3$ is therefore adjacent to at most one of $w_1$ and $w_2$.

Let $y_3$ be a neighbor of $w_3$ in $D_1$. As observed earlier, every vertex in $D_2$ has at most two neighbors in $D_1$. In particular, $w_1$ is not adjacent to $y_3$. If $w_3$ is not adjacent to $w_1$, then $G[\{x_1,x_2,y_1,y_3,w_1, w_3\}]=F_1$, which is a contradiction. Thus, $w_3$ is adjacent to $w_1$, implying that $w_3$ is not adjacent to $w_2$. If $w_2$ is not adjacent to $y_3$, then $G[\{x_1,x_2,y_2,y_3,w_2, w_3\}]=F_1$, which is a contradiction. Hence, $w_2$ is adjacent to $y_3$. If $w_1$ and $w_2$ are not adjacent, then $G[\{x_1,y_1,y_2,y_3,w_1, w_2\}]=P_2 \cart P_3$, which is a contradiction. Hence, $w_1$ and $w_2$ are adjacent, implying that $G[\{x_1,y_2,w_1,w_2,w_3\}]=B$, which is a contradiction. Therefore, $w_2$ is not adjacent to $y_2$.~\smallqed

\medskip

\noindent Recall that $w_2$ is not adjacent to $y_1$. By Claim~\ref{claim5}.1, the vertex $w_2$ is not adjacent to $y_2$. Let $y_4$ be a neighbor of $w_2$ in $D_1$. Since every vertex in $D_2$ has at most two neighbors in $D_1$, the vertex $w_1$ is not adjacent to $y_4$. If $w_1$ is not adjacent to $w_2$, then $G[\{x_1,x_2,y_1,y_4,w_1,w_2\}]=F_1$, which is a contradiction. Hence, $w_1$ is adjacent to $w_2$. 
As $\{x_1, w_1\}$ is not a dominating set of $G$, there is a vertex $w_4$ in $D_2$ that is not adjacent to $w_1$.
Since $G$ is $F_1$-free, the vertex $w_4$ is adjacent to $y_1$ or to $y_2$ or to both $y_1$ and $y_2$.
If $w_4$ is adjacent to $y_1$ and $y_2$, then $G[\{x_1,y_1,y_2,w_1, w_4\}]=K_{2,3}$, which is a contradiction.
Hence, by symmetry, we may assume that $w_4$ is adjacent to $y_1$, but not to $y_2$. 
If $w_4$ is adjacent to $w_2$, then we get the contradiction $G[\{x_1,y_1,y_2,w_1,w_2,w_4\}]=P_2 \cart P_3$.
If $w_4$ is not adjacent to $w_2$, then $G[\{x_1,x_2, y_1,y_4,w_2, w_4\}]=F_1$, which is a contradiction, and completes the proof of Claim~\ref{claim5}.
$\Box$

\medskip

\noindent Let $D_1=\{y_1,\ldots, y_k\}$, and, for $i\in [k]$, let $w_i$ be a neighbor of $y_i$ in $D_2$. 
If $k = 1$, then $\{x_1,y_1\}$ is a dominating set of $G$, which is a contradiction. 
Hence, $k \ge 2$. By Claim~\ref{claim5}, the vertex $y_i$ is the only neighbor of $w_i$ in $D_1$ for $i \in [k]$.
Since $G$ is $F_1$-free, the vertices $w_1,\ldots, w_k$ induce a clique in $G$. 
Thus, by Claim \ref{claim4}, we obtain $k \leq 2$. 
This implies $k=2$.
Since $G$ is $F_1$-free, each neighbor of $y_i$ in $D_2$ is adjacent to every neighbor of $y_{3-i}$ in $D_2$ for $i \in [2]$. Suppose that the vertex $y_1$ has two neighbors $w_1$ and $w_1'$ in $D_2$. Thus, both $w_1$ and $w_1'$ are adjacent to $w_2$. Since $G$ is $D$-free, $w_1$ and $w_1'$ are not adjacent. Since $\{x_1,w_2\}$ is not a dominating set of $G$, the vertex $y_2$ has a neighbor $w_2'$ in $D_2$ that is different from $w_2$ and not adjacent to $w_2$. Thus, $G[\{ w_1,w_1', w_2, w_2', y_1\}]= K_{2,3}$, which is a contradiction, and completes the proof of Theorem~\ref{theorem1}.
$\Box$

\medskip

\noindent We close with a number of conjectures.

\begin{conjecture}\label{conjecture1}
There is a finite set ${\cal F}$ of graphs such that some graph $G$ 
satisfies $\gamma(H) = \gamma_e(H)$ for every induced subgraph $H$ of $G$
if and only if $G$ is ${\cal F}$-free.
\end{conjecture}

\begin{conjecture}\label{conjecture2}
The set ${\cal F}$ in Conjecture \ref{conjecture1} can be chosen such that 
$\gamma(F)=3$ and $\gamma_e(F)=2$ for every graph $F$ in ${\cal F}$.
\end{conjecture}
Similar to the definition of an exponential dominating set, 
Dankelmann et al.~\cite{ddems} define a set $D$ of vertices of a graph $G$ 
to be a {\it porous exponential dominating set} of $G$ if
$w^*_{(G,D)}(u)\geq 1$ for every vertex $u$ of $G$, 
where
$w^*_{(G,D)}(u)=\sum\limits_{v\in D}\left(\frac{1}{2}\right)^{{\rm dist}_G(u,v)-1}$.
They define the {\it porous exponential domination number $\gamma^*_e(G)$} of $G$ 
as the minimum order of a porous exponential dominating set of $G$.
\begin{conjecture}\label{conjecture3}
A graph $G$ satisfies
$\gamma(H)=\gamma_e(H)$ for every induced subgraph $H$ of $G$
if and only if
$\gamma(H)=\gamma_e^*(H)$ for every induced subgraph $H$ of $G$.
\end{conjecture}

\end{document}